\documentclass[11pt, twoside]{article}
\usepackage{fancyhdr}
\usepackage{amssymb}
\usepackage{theorem}
\usepackage{array}
\usepackage{amsmath}
\usepackage{graphicx}
\usepackage{multicol}
\usepackage{comment}

\fontfamily{ptm}\selectfont

\newtheorem{theorem}{Theorem}[section]

\newtheorem{lemma}[theorem]{Lemma}

\newtheorem{definition}[theorem]{Definition}
\newtheorem{remark}[theorem]{Remark}

\newlength{\blspace} 
\newsavebox{\qedbox} 

\sbox{\qedbox}{%
\begin{picture}(5.1,5)%
\put(0,0){\framebox(5,5){}}%
\end{picture}} %
\newcommand{\qed}{\hfill$\blacksquare$}%
\newenvironment{proof}[1][]%
{\begin{trivlist}\item{\textbf{Proof#1.}\hspace{\blspace}}}%
{\qed\end{trivlist}}

\fancypagestyle{plain}{%
   \fancyhead{} 
}

\newcommand{\gln}[1]{\mathrm{GL}_{#1}\mathbb{C}}

\newcommand{\lgln}[1]{\mathfrak{gl}_{#1}\mathbb{C}}

\newcommand{\spinC}{\,\,\,{\makebox[0pt]{$\partial$}\raisebox{1pt}{\makebox[0pt]{$\slash$}}\hspace{3pt}}_\mathbb{C}}
\newcommand{\spinSig}{\,\,\,{\makebox[0pt]{$\partial$}\raisebox{1pt}{\makebox[0pt]{$\slash$}}\hspace{3pt}}_\mathrm{sig}}

\begin{document}

\thispagestyle{plain}

\begin{center}
{\Large \textbf{Signature quantization, representations of compact Lie groups,\\[1mm] and a $q$-analogue of the Kostant partition function}}
\bigskip

Victor \textsc{Guillemin} and Etienne \textsc{Rassart}\\[2mm] 
\textit{Department of Mathematics, Massachusetts Institute of Technology}\\[2mm]
\texttt{\{vwg,rassart\}@math.mit.edu}
\bigskip

\textrm{May 19, 2004}
\end{center}
\bigskip

\begin{abstract}
We discuss some applications of signature quantization to the representation theory of compact Lie groups. In particular, we prove signature analogues of the Kostant formula for weight multiplicities and the Steinberg formula for tensor product multiplicities. Using symmetric functions, we also find, for type $A$, analogues of the Weyl branching rule and the Gelfand-Tsetlin theorem. These analogues involve a $q$-analogue of the Kostant partition function. We show that in type $A$, this $q$-analogue is polynomial in the relative interior of the cells of a complex of cones. This chamber complex can be taken to be the same as the chamber complex of the usual Kostant partition function. We present the case of $A_2$ as a detailed example.
\end{abstract}

\section{Introduction}

The results described in this note are closely related to an article of Guillemin, Sternberg and Weitsman \cite{GSW} on signature quantization.

Denoting by $V_\lambda$ the irreducible representations of a complex semisimple Lie group, the work of Guillemin, Sternberg and Weitsman on quantization with respect to the signature Dirac operator involves the ``twisted'' representations $\widetilde{V}_\lambda = V_{\lambda-\delta}\otimes V_\delta$ for strictly dominant weights $\lambda$. They give a formula for the multiplicities of weights in those representations which is very similar to the Kostant multiplicity formula, but involves the $q=2$ specialization of a $q$-analogue $K_q$ of the Kostant partition function, rather than the usual ($q=1$) partition function. This $q$-analogue arises from the work of Agapito \cite{Agapito} and Guillemin, Sternberg and Weitsman \cite{GSW} in the study of the twisted signature of coadjoint orbits

We explore further the structure of these representation. We explain how they decompose into irreducible representations and show that it is possible to decompose a tensor product of twisted representations into twisted representations again. There is a formula very analogous to that of Steinberg for the multiplicities of the factors in the product, which again involves $K_2$. An interesting feature of this formula and the analogue of the Kostant multiplicity formula of \cite{GSW} is that they are free of the $\delta$ factors of the usual formulas for the irreducible representations. For type $A$, we can write down the characters of the $\widetilde{V}_\lambda$ in terms of Schur functions, and we find a branching rule for restricting the representation $\widetilde{V}_\lambda$ of $\gln{k}$ to $\gln{k-1}$. By iterating this rule, we develop a Gelfand-Tsetlin theory for the twisted representations of $\gln{k}$.

Finally, we describe the structure of the $q$-analogue $K_q$ of the Kostant partition function. We show that for the root system $A_n$, this $q$-analogue is polynomial in the relative interior of the cells of a complex of cones, of degree ${n\choose 2}$ with coefficients in $\mathbb{Q}[q]$ of degree ${n+1\choose 2}$. This chamber complex can be taken to be the same as the chamber complex of the usual Kostant partition function. We present the case of $A_2$ as a detailed example.

\subsection{Quantization with respect to the signature Dirac operator}

A symplectic manifold $(M,\omega)$ is \emph{pre-quantizable} if the cohomology class of $\omega$ is an integral class, i.e. is in the image of the map $H^2(M,\mathbb{Z})\rightarrow H^2(M,\mathbb{R})$. This assumption implies the existence of a \emph{pre-quantum structure} on $M$: a line bundle, $\mathbb{L}$, and a connection, $\nabla$, such that $\mathrm{curv}(\nabla)=\omega$. If $g$ is a Riemannian metric compatible with $\omega$, then, from $g$ and $\omega$, one gets an elliptic operation $\spinC\,:\ S^{+} \rightarrow S^{-}$, the \emph{spin-$\mathbb{C}$ Dirac operator}, and, by twisting this operator with $\mathbb{L}$, an operator $\spinC^{\mathbb{L}}\,:\ S^{+}\otimes\mathbb{L} \rightarrow S^{-}\otimes\mathbb{L}$. If $M$ is compact one can ``quantize'' it by associating with it the virtual vector space
\begin{equation}
\label{eqn:Index}
Q(M) = \mathrm{Index}\spinC^{\mathbb{L}}\,.
\end{equation}
Moreover if $G$ is a compact Lie group and $\tau$ a Hamiltonian action of $G$ on $M$ one gets from $\tau$ a representation of $G$ on $Q(M)$ which is well-defined up to isomorphism (independent of the choice of $g$).

The results described in this note are closely related to two theorems in the article \cite{GSW}. In this article the authors study the signature analogue of spin-$\mathbb{C}$ quantization: i.e. they define the virtual vector space~\eqref{eqn:Index} by replacing $\spinC$ by the signature operator $\spinSig$, and prove signature versions of a number of standard theorems about quantized symplectic manifolds. The two theorems we'll be concerned with in this paper are the following.
\begin{enumerate}
\item Let $G=\left(S^1\right)^n$ and let $M$ be a $2n$-dimensional toric variety with moment polytope $\Delta\subseteq\mathbb{R}^n$. Then, for spin-$\mathbb{C}$ quantization, the weights of the representation of $G$ on $Q(M)$ are the lattice points, $\beta\in\Delta\cap\mathbb{Z}^n$, and each weight occurs with multiplicity $1$. For signature quantization the weights are the same; however, the weight $\beta$ occurs with multiplicity $2^n$ if $\beta$ lies in $\mathrm{Int}(\Delta)$, with multiplicity $2^{n-1}$ if it lies on a facet, and, in general, with multiplicity $2^{n-i}$ if it lies on $i$ facets. Further details can be found in the work of Agapito \cite{Agapito}.
\item Let $G$ be a compact simply connected Lie group, $\lambda$ a dominant weight and $O_\lambda=M$ the coadjoint orbit of $G$ through $\lambda$. In the spin-$\mathbb{C}$ theory, the representation of $G$ on $Q(M)$ is the unique irreducible representation $V_\lambda$ of $G$ with highest weight $\lambda$; however, in the signature theory, it is the representation
\begin{equation}
\widetilde{V}_\lambda = V_{\lambda-\rho}\otimes V_\rho\,,
\end{equation}
where $\rho$ is half the sum of the positive roots. (This is modulo the proviso that $\lambda-\rho$ be dominant.)
\end{enumerate}

The article \cite{GSW} also contains a signature version of the Kostant multiplicity formula. We recall that the Kostant multiplicity formula computes the multiplicity with which a weight, $\mu$, of $T$ occurs in $V_\lambda$ by the formula
\begin{equation}
\label{eqn:KMF}
\sum_{\sigma\in\mathcal{W}}(-1)^{|\sigma|}\,K(\sigma(\lambda+\rho)-(\mu+\rho))
\end{equation}
where $\mathcal{W}$ is the Weyl group, $|\sigma|$ is the length of $\sigma$ in $\mathcal{W}$, and $K$, the \emph{Kostant partition function} (described below in Definition~\ref{def:KPF}). The signature version of the Kostant multiplicity formula computes the multiplicity $\widetilde{m}_\lambda(\mu)$ with which the weight $\mu$ appears in $\widetilde{V}_\lambda$ by a similar formula:
\begin{equation}
\widetilde{m}_\lambda(\mu) = \sum_{\sigma\in\mathcal{W}}(-1)^{|\sigma|}\,K_2(\sigma(\lambda)-\mu)
\end{equation} 
where $K_2$ is the $q=2$ specialization of a new $q$-analogue of the Kostant partition function, described below.

Our initial goal in writing this paper was to give a purely algebraic derivation of this result; however we noticed that there are $\widetilde{V}_\lambda$ analogues of a number of other basic formulas in the representation theory of compact semisimple Lie groups, in particular, an analogue of the Steinberg formula and, for $\gln{k}$, analogues of the Weyl branching rule and the Gelfand-Tsetlin theorem. \nocite{RassartThesis}

\subsection{The Kostant partition function and its $q$-analogues}

We start by introducing the Kostant partition function.

\begin{definition} 
\label{def:KPF} 
The \emph{Kostant partition function} for a root system $\Phi$, given a choice of positive roots $\Phi_+$, is the function 
\begin{equation}
K(\mu) = \Big|\Big\{{(k_\alpha)}_{\alpha\in\Phi_+}\in\mathbb{N}^{|\Phi_+|}\ : \ \sum_{\alpha\in\Phi_+}k_\alpha\alpha = \mu\Big\}\Big|\,,
\end{equation}
i.e. $K(\mu)$ is the number of ways that $\mu$ can be written as a sum of positive roots (see \cite{FultonHarris}).
\end{definition}

Note that $K(\mu)$ can also be computed as the number of integer points inside the polytope
\begin{equation}
Q_\mu = \Big\{{(k_\alpha)}_{\alpha\in\Phi_+}\in\mathbb{R}_{\geq 0}^{|\Phi_+|}\ : \ \sum_{\alpha\in\Phi_+}k_\alpha\alpha = \mu\Big\}\,. 
\end{equation} 

We can write down a generating function for the $K(\mu)$ that is very similar to Euler's generating function for the number of partitions (see \cite[Section 25.2]{FultonHarris}):
\begin{equation}
\sum_{\mu}K(\mu)e^\mu = \prod_{\alpha\in\Phi_+}\frac{1}{1-e^\alpha}\,.
\end{equation}

The classical $q$-analogue $\widehat{K}_q(\mu)$ of $K(\mu)$, due to Lusztig \cite{Lusztig}, keeps track of how many times the roots appear:
\begin{equation}
\widehat{K}_q(\mu) = \sum_{{(k_\alpha)}_\alpha\in Q_\mu}q^{\sum k_\alpha}\,,
\end{equation}
corresponding to the generating function
\begin{equation}
\sum_{\mu}\widehat{K}_q(\mu)e^\mu = \prod_{\alpha\in\Phi_+}\left(\sum_{m\geq 0}q^me^{m\alpha}\right) = \prod_{\alpha\in\Phi_+}\frac{1}{1-qe^\alpha} \,.
\end{equation}


The $q$-analogue $K_q(\mu)$ that interests us here is the one that counts the integer points of $Q_\mu$ according to how many of the $k_\alpha$'s are nonzero:
\begin{equation}
K_q(\mu) = \sum_{{(k_\alpha)}_\alpha\in Q_\mu}q^{|\{k_\alpha > 0\}|}\,.
\end{equation}

In terms of generating functions, this translates to
\begin{equation}
\sum_{\mu}K_q(\mu)e^\mu =  \prod_{\alpha\in\Phi_+}\left(1+q\sum_{m\geq 1}e^{\alpha}\right) = \prod_{\alpha\in\Phi_+}\frac{1+(q-1)e^\alpha}{1-e^\alpha}\,.
\end{equation}


\section{An analogue of the Kostant multiplicity formula}

We are working in the context of a complex semisimple Lie algebra $\mathfrak{g}$ with root system $\Phi$, choice of positive roots $\Phi_+$\,, and Weyl group $\mathcal{W}$\,; $\rho$ is half the sum of the positive roots (or the sum of the fundamental weights). For a dominant weight $\lambda$, we denote by $V_\lambda$ the irreducible representation of $\mathfrak{g}$ with highest weight $\lambda$. We will call a weight $\lambda$ \emph{strictly dominant} if $\lambda-\rho$ is dominant. We will use the notation $\Lambda^{\!+}$ for the set of dominant weights, and $\Lambda^{\!+}_S$ for the set of strictly dominant weights. For a strictly dominant weight, we define the representation
\begin{equation}
\widetilde{V}_\lambda = V_{\lambda-\rho}\otimes V_\rho
\end{equation}
and its character
\begin{equation}
\widetilde{\chi}_{\displaystyle {}_{\lambda}} = \chi_{\displaystyle {}_{V_{\lambda-\rho}\otimes V_\rho}} = \chi_{\displaystyle {}_{\lambda-\rho}}\cdot\chi_{\rho}\,.
\end{equation}

The following theorem of Guillemin, Sternberg, and Weitsman \cite{GSW} provides a formula for the multiplicities of the weights in the weight space decomposition of $\widetilde{V}_\lambda$. This formula is very similar to the Kostant multiplicity formula~\eqref{eqn:KMF}, but uses the $q=2$ specialization of the $q$-analogue of the Kostant partition function $K_q(\mu)$ introduced above, instead of the usual Kostant partition function. The formula for the $\widetilde{V}_\lambda$ multiplicities further distinguishes itself from the Kostant formula by being free of the $\rho$ factors.

\begin{theorem}{\rm\textbf{(Guillemin-Sternberg-Weitsman \cite{GSW})}}\label{thm:GSW}
\ Let $\lambda$ be a strictly dominant weight. Then the multiplicity of the weight $\nu$ in the tensor product $\widetilde{V}_\lambda = V_{\lambda-\rho}\otimes V_{\rho}$ is given by
\begin{equation}
\label{eqn:TwistedMultiplicities}
\widetilde{m}_{\lambda}(\nu) = \mathrm{dim}\,{\big(\widetilde{V}_\lambda\big)}_\nu = \sum_{\omega\in\mathcal{W}}(-1)^{|\omega|}K_2(\omega(\lambda)-\nu)\,,
\end{equation}
where $|\omega|$ is the length of $\omega$ in the Weyl group. 
\end{theorem}

\begin{proof}
We give a simple proof here using the Weyl character formula. This formula expresses the character $\chi_{\displaystyle{}_\lambda}$ of $V_\lambda$ as the quotient
\begin{equation}
\label{eqn:WeylCharacterFormula}
\chi_{\displaystyle{}_\lambda} = \frac{A_{\lambda+\rho}}{A_\rho}\,,
\end{equation}
where $A_\mu = \displaystyle\sum_{\omega\in\mathcal{W}}(-1)^{|\omega|}e^{\omega(\mu)}$\,. For $\rho$, we get the nice expression \cite[Lemma 24.3]{FultonHarris}
\begin{equation}
A_\rho \ =\  \prod_{\alpha\in\Phi_+}\left(e^{\alpha/2}-e^{-\alpha/2}\right) \ =\  e^\rho\prod_{\alpha\in\Phi_+}\left(1-e^{-\alpha}\right)\,,
\end{equation}
which means, in particular, that we get
\begin{equation}
\chi_{\displaystyle{}_\rho} \ =\  \frac{A_{2\rho}}{A_\rho} \ =\  \frac{e^{2\rho}\displaystyle\prod_{\alpha\in\Phi_+}\left(1-e^{-2\alpha}\right)}{e^\rho\displaystyle\prod_{\alpha\in\Phi_+}\left(1-e^{-\alpha}\right)} \ =\  e^\rho\prod_{\alpha\in\Phi_+}\left(1+e^{-\alpha}\right)\,.
\end{equation}

Thus, for $\lambda$ strictly dominant,
\begin{eqnarray}
\widetilde{\chi}_{\displaystyle {}_{\lambda}} \quad=\quad \chi_{\displaystyle {}_{\lambda-\rho}}\cdot\chi_{\displaystyle {}_\rho} & = & \sum_{\omega\in\mathcal{W}}(-1)^{|\omega|}\,e^{\omega(\lambda)}\prod_{\alpha\in\Phi_+}\frac{1+e^{-\alpha}}{1-e^{-\alpha}} \label{eqn:ProductCharacter}\\[3mm]
& = & \sum_{\omega\in\mathcal{W}}(-1)^{|\omega|}\,e^{\omega(\lambda)}\sum_{\mu}K_2(\mu)\,e^{-\mu}\nonumber \\[3mm]
& = & \sum_{\mu}\sum_{\omega\in\mathcal{W}}(-1)^{|\omega|}\,K_2(\mu)\,e^{\omega(\lambda)-\mu}\,.
\end{eqnarray}
Extracting the coefficient of $e^\nu$ on both sides gives~\eqref{eqn:TwistedMultiplicities}.
\end{proof}

The next step will be to use a formula due to Atiyah and Bott for the characters of the $V_\lambda$ and $\widetilde{V}_\lambda$ to break down $\widetilde{V}_\lambda$ into its irreducible components and find their multiplicities. The Atiyah-Bott formula \cite{AtiyahBott1,AtiyahBott2} gives the character of $V_\mu$ as
\begin{equation}
\chi_{\displaystyle {}_{\mu}} = \sum_{\omega\in\mathcal{W}}e^{\omega(\mu)}\prod_{\alpha\in\Phi_+}\frac{1}{1-e^{-\omega(\alpha)}}\,.
\end{equation}

\begin{remark}
We can deduce this formula from the Weyl character formula (equation~\eqref{eqn:WeylCharacterFormula}) by first observing that
\begin{eqnarray}
\prod_{\alpha\in\Delta_+}\left(1-e^{-\omega(\alpha)}\right) & = & \prod_{\begin{scriptsize}\begin{array}{@{}c@{}}\alpha\in\Delta_+ \\ \omega(\alpha)\in\Delta_+\end{array}\end{scriptsize}}\hspace{-3mm}\left(1-e^{-\alpha}\right)\prod_{\begin{scriptsize}\begin{array}{@{}c@{}}\alpha\in\Delta_+ \\ \omega(\alpha)\in\Delta_-\end{array}\end{scriptsize}}\hspace{-3mm}\left(1-e^{\alpha}\right)\nonumber\\[3mm]
& = & \prod_{\begin{scriptsize}\begin{array}{@{}c@{}}\alpha\in\Delta_+ \\ \omega(\alpha)\in\Delta_+\end{array}\end{scriptsize}}\hspace{-3mm}\left(1-e^{-\alpha}\right)\prod_{\begin{scriptsize}\begin{array}{@{}c@{}}\alpha\in\Delta_+ \\ \omega(\alpha)\in\Delta_-\end{array}\end{scriptsize}}\hspace{-3mm}\left(\left(e^{-\alpha}-1\right)e^{\alpha}\right)\nonumber\\[3mm]
& = & (-1)^{|\{\alpha\in\Delta_+\ : \ \omega(\alpha)\in\Delta_-\}|}\prod_{\begin{scriptsize}\begin{array}{@{}c@{}}\alpha\in\Delta_+ \\ \omega(\alpha)\in\Delta_-\end{array}\end{scriptsize}}\hspace{-3mm}e^{\alpha}\prod_{\begin{scriptsize}\begin{array}{@{}c@{}}\alpha\in\Delta_+ \\ \omega(\alpha)\in\Delta_+\end{array}\end{scriptsize}}\hspace{-3mm}\left(1-e^{-\alpha}\right)\prod_{\begin{scriptsize}\begin{array}{@{}c@{}}\alpha\in\Delta_+ \\ \omega(\alpha)\in\Delta_-\end{array}\end{scriptsize}}\hspace{-3mm}\left(1-e^{-\alpha}\right)\nonumber\\[3mm]
& = & (-1)^{|\omega|}e^{\sum\{\alpha\in\Delta_+\ :\ \omega(\alpha)\in\Delta_-\}}\prod_{\alpha\in\Delta_+}\left(1-e^{-\alpha}\right) \label{eqn:AB_W1}
\end{eqnarray}
since the number of positive roots that are sent to negative roots under $\omega$ is the same as the length $|\omega|$ of $\omega$ in the Weyl group.

On the other hand,
\begin{eqnarray}
\delta - \omega(\delta) & = & \frac{1}{2}\sum_{\alpha\in\Delta_+}\alpha \quad-\quad \frac{1}{2}\sum_{\alpha\in\Delta_+}\omega(\alpha)\nonumber\\[4mm]
& = & \frac{1}{2}\sum_{\begin{scriptsize}\begin{array}{@{}c@{}}\alpha\in\Delta_+ \\ \omega(\alpha)\in\Delta_+\end{array}\end{scriptsize}}\hspace{-3mm}\alpha \quad+\quad \frac{1}{2}\sum_{\begin{scriptsize}\begin{array}{@{}c@{}}\alpha\in\Delta_+ \\ \omega(\alpha)\in\Delta_-\end{array}\end{scriptsize}}\hspace{-3mm}\alpha \quad-\quad \left(\frac{1}{2}\sum_{\begin{scriptsize}\begin{array}{@{}c@{}}\alpha\in\Delta_+ \\ \omega(\alpha)\in\Delta_+\end{array}\end{scriptsize}}\hspace{-3mm}\alpha \quad-\quad \frac{1}{2}\sum_{\begin{scriptsize}\begin{array}{@{}c@{}}\alpha\in\Delta_+ \\ \omega(\alpha)\in\Delta_-\end{array}\end{scriptsize}}\hspace{-3mm}\alpha\right)\nonumber\\[4mm]
& = & \sum\{\alpha\in\Delta_+\ : \ \omega(\alpha)\in\Delta_-\}\,. \label{eqn:AB_W2}
\end{eqnarray}

Combining~\eqref{eqn:AB_W1} with~\eqref{eqn:AB_W2} gives
\begin{equation}
\prod_{\alpha\in\Delta_+}\left(1-e^{-\omega(\alpha)}\right) =  (-1)^{|\omega|}e^{\delta-\omega(\delta)}\prod_{\alpha\in\Delta_+}\left(1-e^{-\alpha}\right)\,,
\end{equation}
and we can translate Weyl's character formula into the Atiyah-Bott formula using this equation.
\end{remark}

For any $\omega\in\mathcal{W}$,
\begin{eqnarray}
\chi_{\displaystyle {}_{\rho}} & = & e^\rho\prod_{\alpha\in\Phi_+}\left(1+e^{-\alpha}\right)\nonumber\\[2mm]
& = & e^{\omega(\rho)}\prod_{\alpha\in\Phi_+}\left(1+e^{-\omega(\alpha)}\right)\,,
\end{eqnarray}
since characters are invariant under the Weyl group action. Using this and the Atiyah-Bott formula, we can write\footnote{Alternatively, we can obtain equation~\eqref{eqn:AtiyahBott} from equation~\eqref{eqn:ProductCharacter} by observing that for $\omega\in\mathcal{W}$,
\begin{displaymath}
\omega\cdot\left(\prod_{\alpha\in\Phi_+}\frac{1+e^{-\alpha}}{1-e^{-\alpha}}\right) = \prod_{\alpha\in\Phi_+}\frac{1+e^{-\omega(\alpha)}}{1-e^{-\omega(\alpha)}} = (-1)^{|\omega|}\prod_{\alpha\in\Phi_+}\frac{1+e^{-\alpha}}{1-e^{-\alpha}}\,.
\end{displaymath}
}
\begin{eqnarray}
\label{eqn:AtiyahBott}
\widetilde{\chi}_{\displaystyle{}_\lambda} \quad=\quad \chi_{\displaystyle {}_{\lambda-\rho}}\cdot\chi_{\displaystyle {}_{\rho}} & = & \sum_{\omega\in\mathcal{W}}e^{\omega(\lambda)}\prod_{\alpha\in\Phi_+}\frac{1+e^{-\omega(\alpha)}}{1-e^{-\omega(\alpha)}}\\[3mm]
& = & \sum_{\omega\in\mathcal{W}}e^{\omega(\lambda)}\prod_{\alpha\in\Phi_+}\frac{1}{1-e^{-\omega(\alpha)}}\sum_{I\subseteq\Phi_+}e^{-\omega(\alpha_I)}\nonumber
\end{eqnarray}
where as before, $\alpha_I = \displaystyle\sum_{\alpha\in I}\alpha$\,. This gives
\begin{equation}
\widetilde{\chi}_{\displaystyle{}_\lambda} = \sum_{I\subseteq\Phi_+}\left(\sum_{\omega\in\mathcal{W}}e^{\omega(\lambda-\alpha_I)}\prod_{\alpha\in\Phi_+}\frac{1}{1-e^{-\omega(\alpha)}}\right)\,.
\end{equation}

Letting, $\lambda_I = \lambda-\alpha_I$, we observe that if $\lambda_I$ is dominant, the Atiyah-Bott formula tells us that
\begin{equation}
\sum_{\omega\in\mathcal{W}}e^{\omega(\lambda-\alpha_I)}\prod_{\alpha\in\Phi_+}\frac{1}{1-e^{-\omega(\alpha)}}
\end{equation}
is the character $\chi_{\displaystyle {}_{\lambda_I}}$ of the irreducible representation $V_{\lambda_I}$, so that
\begin{equation}
\widetilde{\chi}_{\displaystyle{}_\lambda} = \sum_{I\subseteq\Phi_+}\chi_{\displaystyle {}_{\lambda_I}} \qquad \textrm{and} \qquad \widetilde{V}_\lambda = V_{\lambda-\rho}\otimes V_{\rho} = \bigoplus_{I\subseteq\Phi_+}V_{\lambda_I}
\end{equation}
if all the $\lambda_I$ are dominant.

Finally, since $\alpha_I$ and $\alpha_{I'}$ can be equal for different subsets $I$ and $I'$, certain highest weights appear multiple times in the above sums. For the weight $\mu = \lambda_I = \lambda-\alpha_I$, we will get $V_\mu$ as many times as we can write $\alpha_I=\lambda-\mu$ as a sum of positive roots, where each positive root appears at most once. Hence
\begin{equation}
\widetilde{V}_\lambda = \sum_{\mu}P(\lambda-\mu)\,V_\mu\,,
\end{equation}
where the sum is over all $\mu$ such that $\mu=\lambda_I$ for some $I$, and $P(\nu)$ is given by
\begin{equation}
\sum_{\nu}P(\nu)e^\nu = \prod_{\alpha\in\Phi_+}\left(1+e^\alpha\right)\,.
\end{equation}

\begin{remark}
David Vogan pointed out to us that this decomposition is well-known and can be deduced from the Steinberg formula. For type $A_n$, the number of distinct $\mu$'s in the above sum is the number of forests of labelled unrooted tree on $n+1$ vertices \cite{KleitmanWinston,Stanley4}. 
\end{remark}

\section{A tensor product formula for the $\widetilde{V}_\lambda$}

We will derive here an analogue of the Steinberg formula for the $\widetilde{V}_\lambda$. Given two representations $\widetilde{V}_\lambda$ and $\widetilde{V}_\mu$, the problem is to determine whether their tensor product $\widetilde{V}_\lambda\otimes\widetilde{V}_\mu$ can be decomposed in terms of $\widetilde{V}_\nu$'s. This is readily seen to be the case, as
\begin{equation}
\widetilde{V}_\lambda\otimes\widetilde{V}_\mu \quad=\quad \left(V_{\lambda-\rho}\otimes V_\rho\right)\otimes\left(V_{\mu-\rho}\otimes V_\rho\right) \quad=\quad \left(V_{\lambda-\rho}\otimes V_\rho\otimes V_{\mu-\rho}\right)\otimes V_\rho\,.
\end{equation}

Breaking up $V_{\lambda-\rho}\otimes V_\rho\otimes V_{\mu-\rho}$ into irreducibles $V_\gamma$ and tensoring each factor with $V_\rho$ yields factors $V_\gamma\otimes V_\rho = \widetilde{V}_{\gamma+\rho}$. Thus for strictly dominant weights $\lambda$ and $\mu$, we can write
\begin{equation}
\widetilde{V}_\lambda\otimes\widetilde{V}_\mu = \sum_{\nu\in\Lambda^{\!+}_S}\widetilde{N}_{\lambda\mu}^{\nu}\widetilde{V}_\nu
\end{equation}
for some nonnegative integers $\widetilde{N}_{\lambda\mu}^{\nu}$.

\begin{theorem} For $\lambda$, $\mu$ and $\nu$ strictly dominant weights, the tensor product multiplicity $\widetilde{N}_{\lambda\mu}^{\nu}$ of $\widetilde{V}_\nu$ in $\widetilde{V}_\lambda\otimes\widetilde{V}_\mu$ is given by 
\begin{equation}
\widetilde{N}_{\lambda\mu}^{\nu} = \sum_{\omega\in\mathcal{W}}\sum_{\sigma\in\mathcal{W}}(-1)^{|\omega\sigma|}\,K_2(\omega(\lambda)+\sigma(\mu)-\nu)\,.
\end{equation}
\end{theorem}

\begin{proof}
Starting from the equation $\widetilde{V}_\lambda\otimes\widetilde{V}_\mu = \displaystyle\sum_{\nu\in\Lambda^{\!+}_S}\widetilde{N}_{\lambda\mu}^{\nu}\widetilde{V}_\nu$, we can use equation~\eqref{eqn:ProductCharacter} to write
\begin{displaymath}
\sum_{\omega\in\mathcal{W}}(-1)^{|\omega|}e^{\omega(\lambda)}\prod_{\alpha\in\Phi_+}\frac{1+e^{-\alpha}}{1-e^{-\alpha}}\ \cdot \ \widetilde{\chi}_{\displaystyle {}_{\mu}}  =  \sum_{\nu\in\Lambda^{\!+}_S}\widetilde{N}_{\lambda\mu}^{\nu}\sum_{\tau\in\mathcal{W}}(-1)^{|\tau|}e^{\tau(\nu)}\prod_{\alpha\in\Phi_+}\frac{1+e^{-\alpha}}{1-e^{-\alpha}}\,.
\end{displaymath}
Cancelling terms and using Theorem~\ref{thm:GSW} to write down the character \raisebox{0.3ex}{$\widetilde{\chi}_{\displaystyle {}_{\mu}}$} yields
\begin{eqnarray*}
\sum_{\omega\in\mathcal{W}}(-1)^{|\omega|}e^{\omega(\lambda)} \ \cdot \ \sum_\beta\sum_{\sigma\in\mathcal{W}}(-1)^{|\sigma|}K_2(\sigma(\mu)-\beta)\,e^{\beta} & = & \sum_{\nu\in\Lambda^{\!+}_S}\widetilde{N}_{\lambda\mu}^{\nu}\sum_{\tau\in\mathcal{W}}(-1)^{|\tau|}e^{\tau(\nu)}\nonumber\\
\sum_{\beta}\sum_{\omega\in\mathcal{W}}\sum_{\sigma\in\mathcal{W}}(-1)^{|\omega|+|\sigma|}\,K_2(\sigma(\mu)-\beta)\,e^{\omega(\lambda)+\beta} & = & \sum_{\nu\in\Lambda^{\!+}_S}\sum_{\tau\in\mathcal{W}}(-1)^{|\tau|}\,\widetilde{N}_{\lambda\mu}^{\nu}\,e^{\tau(\nu)} 
\end{eqnarray*}
Substituting $\gamma=\omega(\lambda)+\beta$ on the left hand side, and $\gamma=\tau(\nu)$ on the right hand side gives
\begin{displaymath}
\sum_{\gamma}\sum_{\omega\in\mathcal{W}}\sum_{\sigma\in\mathcal{W}}(-1)^{|\omega\sigma|}\,K_2(\sigma(\mu)+\omega(\lambda)-\gamma)\,e^{\gamma} = \sum_{\begin{scriptsize}\begin{array}{@{}c@{}}\textrm{$\gamma$ conjugate}\\ \textrm{to a strictly}\\ \textrm{dominant weight}\end{array}\end{scriptsize}}\sum_{\tau\in\mathcal{W}}(-1)^{|\tau|}\,\widetilde{N}_{\lambda\mu}^{\tau^{-1}(\gamma)}\,e^{\gamma}\,, 
\end{displaymath}
and extracting the coefficient of $e^\gamma$ on both sides yields
\begin{equation}
\sum_{\omega\in\mathcal{W}}\sum_{\sigma\in\mathcal{W}}(-1)^{|\omega\sigma|}\,K_2(\sigma(\mu)+\omega(\lambda)-\gamma) = \sum_{\tau\in\mathcal{W}}(-1)^{|\tau|}\,\widetilde{N}_{\lambda\mu}^{\tau^{-1}(\gamma)}\,.
\end{equation}

Now, since $\widetilde{N}_{\lambda\mu}^{\tau^{-1}(\gamma)}$ vanishes unless $\tau^{-1}(\gamma)$ is strictly dominant, all the terms in the sum on the right hand side vanish except for the one where $\tau$ is the identity (i.e. the term where $\gamma=\nu$), and we get the result.
\end{proof}

If we denote by $N_{\lambda\mu}^{\nu}$ the multiplicities of the irreducible representations $V_\nu$ in the tensor product $V_\lambda\otimes V_\mu$, defined by
\begin{equation}
V_\lambda\otimes V_\mu = \sum_{\nu\in\Lambda^{\!+}}N_{\lambda\mu}^{\nu}\,V_\nu\,,
\end{equation}
then we can write down the tensor product multiplicities $\widetilde{N}_{\lambda\mu}^{\nu}$ for the decomposition of $\widetilde{V}_\lambda\otimes\widetilde{V}_\mu$ into $\widetilde{V}_\nu$'s in terms of the $N_{\lambda\mu}^{\nu}$ as follows:
\begin{eqnarray*}
\widetilde{V}_\lambda\otimes\widetilde{V}_\mu & = & V_{\lambda-\rho}\otimes V_\rho\otimes V_{\mu-\rho}\otimes V_\rho\nonumber\\
& = & \left(\left(\sum_{\beta\in\Lambda^{\!+}}N_{\lambda-\rho,\rho}^{\beta}\,V_\beta\right)\otimes V_{\mu-\rho}\right) \otimes V_\rho\nonumber\\
& = & \left(\sum_{\beta\in\Lambda^{\!+}}\sum_{\gamma\in\Lambda^{\!+}}N_{\lambda-\rho,\rho}^{\beta}\,N_{\beta,\mu-\rho}^{\gamma}\,V_\gamma\right) \otimes V_{\rho}\nonumber\\
& = & \sum_{\beta\in\Lambda^{\!+}}\sum_{\gamma\in\Lambda^{\!+}}N_{\lambda-\rho,\rho}^{\beta}\,N_{\beta,\mu-\rho}^{\gamma}\,\widetilde{V}_{\gamma+\rho}\nonumber\\
& = & \sum_{\nu\in\Lambda^{\!+}_S}\sum_{\beta\in\Lambda^{\!+}}N_{\lambda-\rho,\rho}^{\beta}\,N_{\beta,\mu-\rho}^{\nu-\rho}\,\widetilde{V}_\nu\,,
\end{eqnarray*}
so that for strictly dominant $\nu$,
\begin{equation}
\widetilde{N}_{\lambda\mu}^{\nu} = \sum_{\beta\in\Lambda^{\!+}}N_{\lambda-\rho,\rho}^{\beta}\,N_{\beta,\mu-\rho}^{\nu-\rho}\,.
\end{equation}

\begin{remark}
In type $A$, there is a combinatorial interpretation for the coefficients $N_{\lambda\mu}^{\nu}$ in terms of shifted Young tableaux: they are given by a shifted analogue of the Littlewood-Richardson rule (see \cite{Stembridge3}).
\end{remark}

\section{Links with symmetric functions in type $A$}

As for the weight multiplicities and Clebsch-Gordan coefficients, there is a link between the character products \raisebox{0.3ex}{$\widetilde{\chi}_{\displaystyle {}_{\lambda}} = \chi_{\displaystyle {}_{\lambda-\delta}}\cdot\chi_{\displaystyle {}_{\delta}}$} and symmetric functions in type $A$, again in terms of Schur functions.

The character of the irreducible polynomial representation $V_\lambda$ of $\gln{k}$, where we now think of $\lambda$ as a partition with $k$ parts (allowing the empty part) is the Schur function $s_\lambda(x_1,\ldots,x_k)$. We will call a partition \emph{strict} if all its parts are distinct (corresponding to a strictly dominant weight). Thus we have that, for $\gln{k}$,
\begin{equation}
\label{eqn:SchurCharacter}
\widetilde{\chi}_{\displaystyle {}_{\lambda}} \quad=\quad \chi_{\displaystyle {}_{\lambda-\delta}}\cdot\chi_{\displaystyle {}_{\delta}} \quad=\quad s_{\lambda-\delta}(x_1,\ldots,x_k)\,s_\delta(x_1,\ldots,x_k)\,,
\end{equation}
for any strict partition $\lambda$. The weight $\delta$ corresponds to the partition $(k-1,k-2,\ldots,1,0)$.

\begin{remark}
We can also write the characters of $\widetilde{V}_\lambda$ in terms of Hall-Littlewood polynomials. Following \cite[III.1 and III.2]{Macdonald}, for partitions of length at most $k$ with empty parts allowed, let
\begin{equation}
v_m(t) = \prod_{i=1}^{m}\frac{1-t^i}{1-t}
\end{equation}
and define
\begin{equation}
v_\lambda(t) = \prod_{i\geq 0}v_{m_i}(t)
\end{equation}
where $m_i$ is the number of parts of $\lambda$ equal to $i$.

The \emph{Hall-Littlewood polynomials} are the symmetric polynomials defined by
\begin{equation}
P_\lambda(x_1,\ldots,x_k;\,t) = \frac{1}{v_\lambda(t)}R_\lambda(x_1,\ldots,x_k;\,t)\,,
\end{equation}
where $R_\lambda$ is the symmetric polynomial
\begin{equation}
R_\lambda(x_1,\ldots,x_k;\,t) = \sum_{\omega\in\mathfrak{S}_k}\omega\cdot\left(x_1^{\lambda_1}\cdots x_k^{\lambda_k}\prod_{i<j}\frac{x_i-tx_j}{x_i-x_j}\right)\,.
\end{equation}

We can rewrite $R_\lambda$ as
\begin{equation}
R_\lambda(x_1,\ldots,x_k;\,t) = \sum_{\omega\in\mathfrak{S}_k}\omega\cdot\left(x_1^{\lambda_1}\cdots x_k^{\lambda_k}\prod_{i<j}\frac{(1-tx_j/x_i)}{(1-x_j/x_i)}\right)\,.
\end{equation}

For a strict partition $\lambda$ with $k$ parts, $v_\lambda(-1)=1$ and then,
\begin{eqnarray}
P_\lambda(\exp(e_1),\ldots,exp(e_k);\,-1) & = & \sum_{\omega\in\mathfrak{S}_k}\omega\cdot\left(\exp(\lambda_1e_1+\cdots+\lambda_ke_k)\prod_{i<j}\frac{(1+\exp(e_j-e_i))}{(1-\exp(e_j-e_i))}\right)\nonumber\\[3mm]
& = & \sum_{\omega\in\mathfrak{S}_k}\omega\cdot\left(e^\lambda\prod_{\alpha\in\Delta_+}\frac{1+e^{-\alpha}}{1-e^{-\alpha}}\right)\nonumber\\[3mm]
& = & \sum_{\omega\in\mathfrak{S}_k}e^{\omega(\lambda)}\prod_{\alpha\in\Delta_+}\frac{1+e^{-\omega(\alpha)}}{1-e^{-\omega(\alpha)}}\nonumber\\[3mm]
& = & \chi_{\displaystyle {}_{\lambda-\delta}}\cdot\chi_{\displaystyle {}_{\delta}}
\end{eqnarray}
from the Atiyah-Bott formula (equation~\eqref{eqn:AtiyahBott}). So the character product \raisebox{0.3ex}{$\chi_{\displaystyle {}_{\lambda-\delta}}\cdot\chi_{\displaystyle {}_{\delta}}$} can be thought of as the $t=-1$ specialization of the Hall-Littlewood polynomial $P_\lambda$.

The results of the following sections can be deduced from this link with Hall-Littlewood polynomials, but we will rather use the Schur function expression~\eqref{eqn:SchurCharacter} for the characters. This makes the proofs a bit more technical but avoids the heavier machinery of Hall-Littlewood polynomials.
\end{remark}

\section{A branching rule for the $\widetilde{V}_\lambda$ in type $A$}

We have seen that the representations $\widetilde{V}_\lambda$ behave somewhat like irreducible representations, in that tensor products of them can be broken down into direct sums of $\widetilde{V}_\nu$'s again, and that the multiplicities in those decompositions as well as in the weight space decomposition are given by formulas very similar to those of Kostant and Steinberg in the irreducible case. The Weyl branching rule (see \cite{FultonHarris} for example) describes how to restrict a representation $V_\lambda$ from $\gln{k}$ to $\gln{k-1}$. This rule can be applied iteratively and prodides a way to index one-dimensional subspaces of $V_\lambda$ by diagrams (Gelfand-Tsetlin diagrams \cite{GelfandTsetlin}) that is compatible with the weight space decomposition. It is natural to ask whether the representations $\widetilde{V}_\lambda$ of $\gln{k}$ are also well-behaved under restriction, or in another words, if there is an analogue of the Weyl branching rule for the $\widetilde{V}_\lambda$ in type $A$.

For two partitions $\mu=(\mu_1,\ldots,\mu_m)$ and $\gamma=(\gamma_1,\ldots,\gamma_{m-1})$, we say that $\gamma$ interlaces $\mu$, and write $\gamma \lhd \mu$, if
\begin{displaymath}
\mu_1 \geq \gamma_1 \geq \mu_2 \geq \gamma_2 \geq \mu_3 \geq \cdots \geq \mu_{m-1} \geq \gamma_{m-1} \geq \mu_m\,.
\end{displaymath}

For two such partitions $\mu$ and $\gamma$ such that $\gamma \lhd \mu$, we define
\begin{equation}
\nabla(\mu,\gamma) = \big|\big\{i\in\{1,2,\ldots,m-1\}\ : \ \mu_i > \gamma_i > \mu_{i+1}\big\}\big|\,.
\end{equation}
In other words, $\nabla(\mu,\gamma)$ is the number of $\gamma_i$ that are wedged strictly between $\mu_i$ and $\mu_{i+1}$.

\begin{theorem}
The decomposition of the restriction of the representation $\widetilde{V}_\lambda$ of $\gln{k}$ to $\gln{k-1}$ into irreducible representations of $\gln{k-1}$ is given by 
\begin{equation}
\mathrm{Res}_{\gln{k-1}}^{\gln{k}}\,\widetilde{V}_\lambda = \bigoplus_{\nu\in\Lambda^{\!+}_S\,:\,\nu \,\lhd\, \lambda}2^{\nabla(\lambda,\nu)}\,\widetilde{V}_\nu\,.
\end{equation} 
\end{theorem}

\begin{proof}
We will argue using characters and the fact that those can be written in terms of Schur functions. 

We saw above (equation~\eqref{eqn:SchurCharacter}) that the character of the representation $\widetilde{V}_\lambda$ of $\gln{k}$ is the product of Schur functions $s_{\lambda-\delta}(x_1,\ldots,x_k)\,s_\delta(x_1,\ldots,x_k)$. We obtain the character of the restriction of $\widetilde{V}_\lambda$ to $\gln{k-1}$ by setting the last variable $x_k$ equal to $1$. Now if we have a Schur function in two sets of variables $x=(x_1,x_2,\ldots)$ and $y=(y_1,y_2,\ldots)$ with the ordering $x_1<x_2<\cdots<y_1<y_2<\cdots$, then we have the identity
\begin{equation}
s_\lambda(x,y) = \sum_{\mu \subseteq \lambda}s_\mu(x)\,s_{\lambda/\mu}(y)\,.
\end{equation}
This is simply saying that we get a semistandard Young tableau of shape $\lambda$ with entries in $x$ and $y$ by first filling a subtableau $\mu$ with entries in $x$ and then the remaining skew-shape with entries from $y$. In our case, with $x=(x_1,\ldots,x_{k-1})$ and $y=1$, this yields
\begin{equation}
s_\lambda(x_1,\ldots,x_{k-1},1) = \sum_{\mu \subseteq \lambda}s_\mu(x_1,\ldots,x_{k-1})\,s_{\lambda/\mu}(1)\,.
\end{equation}
However, we have $s_{\lambda/\mu}(1) = $ unless $\lambda/\mu$ is a horizontal strip, in which case it is equal to $1$, and also $s_\mu(x_1,\ldots,x_{k-1}) = 0$\, if $\mu$ has $k$ parts or more. Hence
\begin{equation}
s_\lambda(x_1,\ldots,x_{k-1},1) = \sum_{\mu}s_\mu(x_1,\ldots,x_{k-1})
\end{equation}  
where the sum is over all $\mu$ with at most $k$ parts that can be obtained from $\lambda$ by removing a horizontal strip. The set of such $\mu$'s is seen to be the set of partitions that interlace $\lambda$, so
\begin{equation}
\label{eqn:SchurTwoVariableSets}
s_\lambda(x_1,\ldots,x_{k-1},1) = \sum_{\mu\,\lhd\,\lambda}s_\mu(x_1,\ldots,x_{k-1})\,.
\end{equation}
Also, in the $x_i=exp(e_i)$ coordinates,
\begin{eqnarray}
\chi_{\displaystyle {}_{\delta}} & = & e^\delta\prod_{\alpha\in\Delta_+}\left(1+e^{-\alpha}\right)\nonumber\\[2mm]
& = & x_1^{k-1}x_2^{k-2}\cdots x_{k-1}\prod_{1\leq i<j\leq k}\left(1+\frac{x_j}{x_i}\right)\nonumber\\[2mm]
& = & \prod_{1\leq i<j\leq k}(x_i+x_j)\,.
\end{eqnarray}
This can also be deduced from the classical definition of the Schur functions in terms of determinants \cite[Section 7.15]{EC2}, since $s_{\delta}(x_1,\ldots,x_k)$ is the ratio between the Vandermonde determinant in variables $x_1^2,\ldots,x_k^2$ and the Vandermonde determinant in $x_1,\ldots,x_k$. Thus,
\begin{equation}
s_{\lambda-\delta}(x_1,\ldots,x_{k-1},1)s_\delta(x_1,\ldots,x_{k-1},1) = \sum_{\mu\,\lhd\,\lambda-\delta}s_\mu(x_1,\ldots,x_{k-1})\hspace{-3mm}\prod_{1\leq i<j\leq k-1}\hspace{-3mm}(x_i+x_j)\prod_{i=1}^{k-1}(x_i+1)\,.
\end{equation}

We recognize the product $\prod_{1\leq i<j\leq k-1}(x_i+x_j)$ as the Schur function $s_\delta(x_1,\ldots,x_{k-1})$ (where $\delta$ now corresponds to the partition $(k-2,k-3,\ldots,1,0)$ with $k-1$ parts), and the product $\prod_{i=1}^{k-1}(x_i+1)$ as the sum $(e_0+e_1+\cdots+e_{k-1})$ of elementary symmetric functions in the variables $x_1,\ldots,x_{k-1}$. A dual version of the Pieri rule \cite[Section 7.15]{EC2} describes how to break down the product of a Schur function with an elementary symmetric function into Schur functions:
\begin{equation}
s_\mu\,e_m = \sum_{\nu}s_\nu\,,
\end{equation}
where the sum is over all $\nu$ obtained from $\mu$ by adding a vertical strip of size $m$, i.e. over the $\nu$ such that $\mu\subseteq\nu$ and the skew-shape $\nu/\mu$ consists of $m$ boxes, no two of which are in the same row. As we are working in ${k-1}$ variables, the $s_\nu$ with more than ${k-1}$ parts vanish, so we can add the further constraint that the vertical strip be confined to the first $k-1$ rows (we will say such a vertical strip has height at most $k-1$). This gives
\begin{eqnarray}
s_{\lambda-\delta}(x_1,\ldots,x_{k-1},1)s_\delta(x_1,\ldots,x_{k-1},1) & = & \sum_{\mu\,\lhd\,\lambda-\delta}\sum_{\nu}s_\nu(x_1,\ldots,x_{k-1})\,s_\delta(x_1,\ldots,x_{k-1})\nonumber\\[2mm]
\widetilde\chi_{\displaystyle {}_{\lambda}}(x_1,\ldots,x_{k-1},1) & = & \sum_{\mu\,\lhd\,\lambda-\delta}\sum_{\nu}\widetilde\chi_{\displaystyle {}_{\nu+\delta}}(x_1,\ldots,x_{k-1})
\end{eqnarray}
where the sum is over all the $\nu$ that can be obtained from $\mu$ by adding a vertical strip of size and height at most $k-1$. We can rewrite this as
\begin{equation}
\widetilde\chi_{\displaystyle {}_{\lambda}}(x_1,\ldots,x_{k-1},1) = \sum_{\mu\,\lhd\,\lambda-\delta}\sum_{\nu}\widetilde\chi_{\displaystyle {}_{\nu}}(x_1,\ldots,x_{k-1})
\end{equation}
where the sum is over all strict partitions $\nu$ such that $\nu-\delta$ can be obtained from $\mu$ by adding a vertical strip of size and height at most $k-1$. Since the $s_\nu s_\delta$ are linearly independent, we can lift this to the level of representations to get
\begin{equation}
\label{eqn:Restriction1}
\mathrm{Res}_{\gln{k-1}}^{\gln{k}}\,\widetilde{V}_\lambda = \bigoplus_{\mu\,\lhd\,\lambda-\delta}\bigoplus_{\nu}\widetilde{V}_\nu\,,
\end{equation}
with the sum over the same set of $\nu$ as before.

In order to compute the multiplicity of a given $\widetilde{V}_\nu$ in $\mathrm{Res}_{\gln{k-1}}^{\gln{k}}\,\widetilde{V}_\lambda$, we define, for strict partitions $\lambda$ and $\nu$, $n(\lambda,\nu)$ to be the number of ways that $\nu-\delta$ can be obtained by adding a vertical strip of size and height at most $k-1$ to some partition $\mu$ such that $\mu\,\lhd\,\lambda-\delta$, so that
\begin{equation}
\label{eqn:Restriction2}
\widetilde{V}_\lambda = \bigoplus_{\nu\in\Lambda^{\!+}_S}n(\lambda,\nu)\,\widetilde{V}_\nu\,.
\end{equation}
Note that $\delta$ has two different meanings here: for the group $\gln{k}$, it corresponds to the partition $(k-1,k-2,\ldots,1,0)$, while for $\gln{k-1}$, it corresponds to the partition $(k-2,k-3,\ldots,1,0)$. To avoid confusion, we will denote the latter by $\delta'$.

The condition $\mu\,\lhd\,\lambda-\delta$ means that
\begin{displaymath}
\lambda_1-(k-1) \geq \mu_1 \geq \lambda_2-(k-2) \geq \mu_2 \geq \cdots \geq \lambda_{k-1}-1 \geq \mu_{k-1} \geq \lambda_k\,.
\end{displaymath}
Replacing $\mu_i$ by $\mu_i+\delta'_i = \mu_i+(k-1-i)=$ gives
\begin{eqnarray*}
\lambda_1-1 \geq & \mu_1+(k-2) & \geq \lambda_2\\
\lambda_2-1 \geq & \mu_2+(k-1) & \geq \lambda_3\\
\ldots & \ldots & \ldots\\
\lambda_{k-1}-1 \geq & \mu_{k-1}+(0) & \geq \lambda_k\,.
\end{eqnarray*}
These equations mean that the $i$-th part of $\mu'=\mu+\delta'$ is at least as large as the $(i+1)$-th part of $\lambda$ and smaller than the $i$-th part of $\lambda$. In other words, the skew-shape $\lambda/\mu'$ is a horizontal strip with a least a box in each row, or equivalently, $\mu'\,\lhd\,\lambda$ with the further constraints $\mu'_i<\lambda_i$ for all $1\leq i\leq k-1$. Adding a vertical strip to $\mu$ to get $\nu-\delta$ is the same as adding a vertical strip to $\mu'$ to get $\nu$, provided that we only allow adding vertical strips to $\mu'$ that result in a strict partition. It is then clear that by adding such a vertical strip to $\mu'$, we get a strict partition $\nu$ such that $\lambda/\nu$ is a horizontal strip. Conversely, it is also clear that for any strict $\nu$ such that $\lambda/\nu$ is a horizontal strip, there is a $\mu'$ such that $\nu$ can be obtained from $\mu'$ by adding a vertical strip. So the only summands $\widetilde{V}_\nu$ for which $n(\lambda,\nu)\neq 0$ in the decomposition~\eqref{eqn:Restriction2} are those for which $\nu\,\lhd\,\lambda$.

Given such a $\nu$, we will compute $n(\lambda,\nu)$ by constructing row by row the strict partitions $\mu'=\mu+\delta'$ from which we can obtain $\nu$. Given $\nu_i$, there are three cases to consider for the possible $\mu_i'$\,:
\begin{itemize}
\itemsep0.5ex
\item $\nu_i=\lambda_i$. In this case, since we must have $\mu_i'<\lambda_i$, it has to be that $\mu_i'=\lambda_i-1$ and that we have a box in row $i$ of the vertical strip. So there is only one choice for $\mu_i'$\,.
\item $\nu_i=\lambda_{i+1}$. Then we must have $\mu_i=\lambda_{i+1}\leq \mu_i'\leq \nu_i$ and therefore $\mu_i'=\nu_i$, so we don't have a box in row $i$ of the vertical strip. Again, there is only one choice for $\mu_i'$ in this case.
\item $\lambda_i>\nu_i>\lambda_{i+1}$. Then we can either have $\mu_i'=\nu_i-1$ and have a box from the vertical strip in row $i$, or have $\mu_i'=\nu_i$ and have no box from the vertical strip in row $i$. So there are two possibilities for $\mu_i'$ in this case.
\end{itemize}
We have to show that any choice of $\mu_i'$ that we make gives rise to a strict partition (by construction, it is clear that $\mu'\,\lhd\,\lambda$). If for some $i$ we had $\mu_i'=\mu_{i+1}'$, then because $\lambda_{i+1}$ is at least $\mu_{i+1}'+1$, this would mean that $\lambda_i$ is at least $\mu_i'+2$, since $\lambda_i>\lambda_{i+1}$. But then $\lambda/\mu'$ contains two boxes in the same column: the box after box $\mu_i'$ in row $i$, and the box after box $\mu_i'=\mu_{i+1}'$ in row $i+1$, which contradicts the fact that $\mu'\,\lhd\,\lambda$ (or equivalently, that $\lambda/\mu'$ is a horizontal strip). Hence we get two choices for each instance of a pattern of the form $\lambda_i>\nu_i>\lambda_{i+1}$. We called the number of such instances above $\nabla(\lambda,\nu)$. Since the choices at each row are independent, we have
\begin{equation}
n(\lambda,\nu) = \left\{\begin{array}{ll}
2^{\nabla(\lambda,\nu)} & \textrm{if $\nu\,\lhd\,\lambda$ and $\nu\in\Lambda^{\!+}_S$\,,}\\[2mm]
0 & \textrm{otherwise.}\end{array}\right.
\end{equation}
from which the proposed expression for the branching rule follows.   
\end{proof}

\section{A Gelfand-Tsetlin theory for the $\widetilde{V}_\lambda$ in type $A$}

After restricting to $\gln{k-1}$, we can further restrict to $\gln{k-2}$. From now on, we will assume that all partitions are strict. We can write
\begin{eqnarray}
\mathrm{Res}_{\gln{k-2}}^{\gln{k}}\,\widetilde{V}_\lambda & = & \mathrm{Res}_{\gln{k-2}}^{\gln{k-1}}\left(\mathrm{Res}_{\gln{k-1}}^{\gln{k}}\,\widetilde{V}_\lambda\right)\nonumber\\
& = & \mathrm{Res}_{\gln{k-2}}^{\gln{k-1}}\,\left(\bigoplus_{\nu \,\lhd\, \lambda}2^{\nabla(\lambda,\nu)}\,\widetilde{V}_\nu\right)\nonumber\\
& = & \bigoplus_{\nu \,\lhd\, \lambda}2^{\nabla(\lambda,\nu)}\,\mathrm{Res}_{\gln{k-2}}^{\gln{k-1}}\,\widetilde{V}_{\nu}\nonumber\\
& = & \bigoplus_{\nu \,\lhd\, \lambda}2^{\nabla(\lambda,\nu)}\left(\bigoplus_{\mu \,\lhd\, \nu}2^{\nabla(\nu,\mu)}\,\widetilde{V}_\mu\right)\nonumber\\
& = & \bigoplus_{\mu \,\lhd\, \nu \,\lhd\, \lambda}2^{\nabla(\lambda,\nu)+\nabla(\nu,\mu)}\,\widetilde{V}_\mu\,.
\end{eqnarray}

Denoting by $\lambda^{(m)}=\lambda^{(m)}_1\geq\cdots\geq\lambda^{(m)}_m\geq0$ the strict partitions indexing the representations $\widetilde{V}$ of $\gln{m}$, we can iterate the branching rule until we get to $\gln{1}$\,:
\begin{equation}
\mathrm{Res}_{\gln{1}}^{\gln{k}}\,\widetilde{V}_{\lambda} = \bigoplus_{\lambda^{(1)}\,\lhd\,\cdots\,\lhd\,\lambda^{(k)}\,=\,\lambda}\hspace{-3mm}2^{\nabla(\lambda^{(k)},\lambda^{(k-1)})+\nabla(\lambda^{(k-1)},\lambda^{(k-2)})+\cdots+\nabla(\lambda^{(2)},\lambda^{(1)})}\,V_{\lambda^{(1)}}\,.
\end{equation}

\newcommand{\dddots}{\raisebox{-1.1mm}[0cm][0cm]{$\cdot$}\,\raisebox{0mm}[0cm][0cm]{$\cdot$}\,\raisebox{1.1mm}[0cm][0cm]{$\cdot$}}

We will call a sequence of strict partitions of the form ${\lambda^{(1)}\,\lhd\,\cdots\,\lhd\,\lambda^{(k)}\,=\,\lambda}\ $ a \emph{twisted Gelfand-Tsetlin diagram} for $\lambda$, which can be viewed schematically as
\begin{equation}
\renewcommand{\arraystretch}{1.5}
\begin{array}{c@{\hspace{2mm}}c@{\hspace{2mm}}c@{\hspace{2mm}}c@{\hspace{2mm}}c@{\hspace{2mm}}c@{\hspace{2mm}}c@{\hspace{2mm}}c@{\hspace{2mm}}c}
\lambda^{(k)}_1 & & \lambda^{(k)}_2 & & \cdots & & \lambda^{(k)}_{k-1} & & \lambda^{(k)}_k \\
& \lambda^{(k-1)}_{1} & & \lambda^{(k-1)}_{2} & & \cdots & & \lambda^{(k-1)}_{k-1} & \\
& & \ddots & & \vdots & & \ \raisebox{0.5mm}{\dddots}  & & \\ 
& & & \lambda^{(2)}_{1} & & \lambda^{(2)}_{2} & & & \\
& & & & \lambda^{(1)}_1 & & & & 
\end{array}
\end{equation}
with $\lambda^{(k)}_j=\lambda_j$ and each $\lambda^{(i)}_j$ is a nonnegative integer satisfying
\begin{equation}
\lambda^{(i)}_j > \lambda^{(i)}_{j+1}
\end{equation}
and
\begin{equation}
\lambda^{(i+1)}_{j}\  \geq \ \lambda^{(i)}_{j} \ \geq \ \lambda^{(i+1)}_{j+1} 
\end{equation}
for all $1\leq j\leq i$, $1\leq i\leq k-1$.

Let $\widetilde{V}_\mathcal{D}$ be the subspace of $\widetilde{V}_{\lambda}$ corresponding to a twisted Gelfand-Tsetlin diagram $\mathcal{D}$. This subspace has dimension $2^{\nabla(\mathcal{D})}$, where
\begin{equation}
\nabla(\mathcal{D}) = \nabla(\lambda^{(k)},\lambda^{(k-1)})+\nabla(\lambda^{(k-1)},\lambda^{(k-2)})+\cdots+\nabla(\lambda^{(2)},\lambda^{(1)})\,.
\end{equation}

We can also think of $\nabla(\mathcal{D})$ as the number of triangles
\begin{displaymath}
\begin{array}{@{}c@{\hspace{2mm}}c@{\hspace{2mm}}c@{}}
\lambda^{(i)}_j & & \lambda^{(i)}_{j+1}\\[3mm]
& \lambda^{(i+1)}_{j} & 
\end{array}
\end{displaymath}
with strict inequalities $\lambda^{(i+1)}_{j} > \lambda^{(i)}_{j} > \lambda^{(i+1)}_{j+1}$ in the diagram $\mathcal{D}$. 

We show here that $\widetilde{V}_\mathcal{D}$ lies completely within the same weight space of the weight space decomposition of $\widetilde{V}_\lambda$.

Consider $\gln{k}$ with its subgroup $T_k$ of invertible diagonal matrices, and also its Lie algebra $\lgln{k}$ and its ``Cartan'' subalgebra $\mathfrak{t}_k$ of diagonal matrices. We have the natural basis in which weights are usually written
\begin{displaymath}
\renewcommand{\arraystretch}{0.7}
J_1=\left(\begin{array}{c@{\hspace{2mm}}c@{\hspace{2mm}}c@{\hspace{2mm}}c@{\hspace{2mm}}c}
1 & & & &\\
 & 0& & &\\
 & & 0& &\\[-1mm]
 & & & \ddots &\\
 & & & & 0
\end{array}\right),\quad
J_2=\left(\begin{array}{c@{\hspace{2mm}}c@{\hspace{2mm}}c@{\hspace{2mm}}c@{\hspace{2mm}}c}
0 & & & &\\
 & 1& & &\\
 & & 0& &\\[-1mm]
 & & & \ddots &\\
 & & & & 0
\end{array}\right), \quad\ldots\quad,\quad 
J_k=\left(\begin{array}{c@{\hspace{2mm}}c@{\hspace{2mm}}c@{\hspace{2mm}}c@{\hspace{2mm}}c}
0 & & & &\\
 & 0& & &\\
 & & 0& &\\[-1mm]
 & & & \ddots &\\
 & & & & 1
\end{array}\right)
\renewcommand{\arraystretch}{1.0}
\end{displaymath}
for $\mathfrak{t}_k$, and also the basis
\begin{displaymath}
\renewcommand{\arraystretch}{0.7}
I_1=\left(\begin{array}{c@{\hspace{2mm}}c@{\hspace{2mm}}c@{\hspace{2mm}}c@{\hspace{2mm}}c}
1 & & & & \\
 & 0& & & \\
 & & 0& & \\[-1mm]
 & & & \ddots &\\
 & & & & 0
\end{array}\right),\quad 
I_2=\left(\begin{array}{c@{\hspace{2mm}}c@{\hspace{2mm}}c@{\hspace{2mm}}c@{\hspace{2mm}}c}
1 & & & &\\
 & 1& & &\\
 & & 0& &\\[-1mm]
 & & & \ddots &\\
 & & & & 0
\end{array}\right), \quad\ldots\quad,\quad 
I_k=\left(\begin{array}{c@{\hspace{2mm}}c@{\hspace{2mm}}c@{\hspace{2mm}}c@{\hspace{2mm}}c}
1 & & & &\\
 & 1& & &\\
 & & 1& &\\[-1mm]
 & & & \ddots &\\
 & & & & 1
\end{array}\right)
\renewcommand{\arraystretch}{1.0}
\end{displaymath}
which is more convenient and which we will use to do the computation. We will simply have to remember that $J_i=I_i-I_{i-1}$ to get the weights in their usual form at the end.

We will think of the groups $\gln{k}$ as included into one another by identifying $\gln{m}$ with
\begin{displaymath}
\renewcommand{\arraystretch}{2.0}
\left(\begin{array}{c|c}
\gln{m} & \mathbf{0}\\
\hline
\mathbf{0} & \mathit{id}_{k-m}
\end{array}\right)
\renewcommand{\arraystretch}{1.0}
\end{displaymath}

Consider the element $I\in\lgln{m}$ and a representation $\widetilde{V}_\mu$ of $\gln{m}$. Then we have the representation $\gln{k} \rightarrow \mathfrak{gl}(V_\mu\otimes V_\rho)$. For $v\in V_{\mu-\rho}$ and $w\in V_\rho$, we have
\begin{eqnarray*}
I\cdot(v\otimes w) & = & (I\cdot v)\otimes w + v\otimes(I\cdot w)\\
& = & \left(\left(\sum_{j=1}^{m}(\mu-\rho)_j\right)v\right)\otimes w + v\otimes\left(\left(\sum_{j=1}^{m}\rho_j\right)w\right)\\
& = & \left(\sum_{j=1}^{m}\left((\mu-\rho)_j+\rho_j\right)\right)\,v\otimes w\\
& = & \left(\sum_{j=1}^{m}\mu_j\right)\,v\otimes w\,, 
\end{eqnarray*}
since $V_{\mu-\rho}$ has highest weight $\mu-\rho$ and $V_\rho$ has highest weight $\rho$. So $I\in\lgln{m}$ gets represented as $(\sum_{j=1}^{m}\mu_j)\,I$ in $\widetilde{V}_\mu$.

In general, for 
\begin{displaymath}
\mathrm{Res}_{\gln{m}}^{\gln{k}}\,\widetilde{V}_\lambda = \bigoplus_{\lambda^{(m)}\,\lhd\,\cdots\,\lhd\,\lambda^{(k)}=\lambda}2^{\nabla(\lambda^{(k)},\lambda^{(k-1)})+\nabla(\lambda^{(k-1)},\lambda^{(k-2)})+\cdots+\nabla(\lambda^{(m+1)},\lambda^{(m)})}\,\widetilde{V}_{\lambda^{(m)}}\,,
\end{displaymath}
we will find that $I\in\lgln{m}$ gets represented as $(\sum^{m}_{i=1}\lambda^{(m)}_i)\,I$ in $\widetilde{V}_{\lambda^{(m)}}$.

Therefore, in the basis $I_1,\,\ldots\,,I_k$, the subspace $\widetilde{V}_\mathcal{D}$ corresponding to a twisted Gelfand-Tsetlin diagram $\mathcal{D}$ has weight
\begin{displaymath}
\Big(\sum_{i=1}^{1}\lambda^{(1)}_i\,, \ \sum_{i=1}^{2}\lambda^{(2)}_i\,, \ \ldots\,, \ \sum_{i=1}^{k}\lambda^{(k)}_i\Big)
\end{displaymath}
or
\begin{equation}
\label{eqn:twistedGTweight0}
\Big(\sum_{i=1}^{1}\lambda^{(1)}_i\,, \ \sum_{i=1}^{2}\lambda^{(2)}_i-\sum_{i=1}^{1}\lambda^{(1)}_i\,, \ \ldots\,, \ \sum_{i=1}^{k}\lambda^{(k)}_i-\sum_{i=1}^{k-1}\lambda^{(k-1)}_i\Big)
\end{equation}
in the usual basis $J_1,\ldots,J_k$.

In other words, $\widetilde{V}_\mathcal{D}\subseteq\big(\widetilde{V}_\lambda\big)_{\beta}$ if 
\begin{equation}
\beta_m = \sum_{i=1}^{m}\lambda^{(m)}_i-\sum_{i=1}^{m-1}\lambda^{(m-1)}_i\,,
\label{eqn:twistedGTweight1}
\end{equation}
or, equivalently,
\begin{equation}
\beta_1 + \cdots + \beta_m = \sum_{i=1}^{m}\lambda^{(m)}_i\,.
\label{eqn:twistedGTweight2}
\end{equation}

Hence twisted Gelfand-Tsetlin diagrams for $\lambda$ correspond to the same weight if all their row sums are the same. So we have proved the following analogue of the Gelfand-Tsetlin theorem \cite{GelfandTsetlin}.

\begin{theorem}
Let $\lambda=(\lambda_1,\ldots,\lambda_k)$ be a strictly dominant weight. The dimension of the representation $\widetilde{V}_\lambda$ of $\gln{k}$ is given by
\begin{equation}
\dim\,\widetilde{V}_\lambda = \sum_{\mathcal{D}}2^{\nabla(\mathcal{D})}
\end{equation} 
where the sum is over all twisted Gelfand-Tsetlin diagrams with top row $\lambda$.

Furthermore, the multiplicity $\widetilde{m}_\lambda(\beta)$ of the weight $\beta$ in $\widetilde{V}_\lambda$ is given by
\begin{equation}
\widetilde{m}_\lambda(\beta) = \dim\,{\big(\widetilde{V}_\lambda\big)}_\beta = \sum_{\mathcal{D}}2^{\nabla(\mathcal{D})}
\end{equation}
where the sum is over all twisted Gelfand-Tsetlin diagrams with top row $\lambda$ and row sums satisfying equation~(\ref{eqn:twistedGTweight1}) (or (\ref{eqn:twistedGTweight2})). 
\end{theorem}

\begin{remark}
We can also prove that $\widetilde{V}_\mathcal{D}$ lies completely within a weight space of $\widetilde{V}_\lambda$ using characters. The Schur function identity~\eqref{eqn:SchurTwoVariableSets}
\begin{equation}
s_\lambda(x,y) = \sum_{\mu\subseteq\lambda}s_\mu(x)\,s_{\lambda/\mu}(y)
\end{equation}
in the two sets of variables $x$ and $y$ gives
\begin{eqnarray}
s_\lambda(x_1,\ldots,x_{k-1},x_k) & = & \sum_{\mu}s_\mu(x_1,\ldots,x_{k-1})\,x_k^{|\lambda/\mu|}\nonumber\\
& = & \sum_{\mu}s_\mu(x_1,\ldots,x_{k-1})\,x_k^{|\lambda|-|\mu|}
\end{eqnarray}
where the sum is over all $\mu$ such that $\lambda/\mu$ is a horizontal strip (or equivalently, $\mu$ such that $\mu\,\lhd\,\lambda$). We also have
\begin{equation}
\prod_{i=1}^{k-1}(x_i+x_k) = x_k^{k-1}e_0 + x_k^{k-2}e_2 + \cdots + x^ke_{k-2} + e_{k-1}\,,
\end{equation}
where as before, the $e_m$ are the elementary symmetric functions in the variables $x_1,\ldots,x_k$. This gives
\begin{eqnarray*}
s_\lambda(x_1,\ldots,x_k) & = & \sum_{\mu\,\lhd\,\lambda-\delta}s_{\mu}(x_1,\ldots,x_{k-1})x_k^{|\lambda|-|\delta|-|\mu|}\,s_\delta(x_1,\ldots,x_{k-1})\sum_{m=0}^{k-1}x_k^{k-1-m}\,e_m \\
& = & \sum_{\mu\,\lhd\,\lambda-\delta}\sum_{\nu}s_\nu(x_1,\ldots,x_{k-1})\,s_\delta(x_1,\ldots,x_{k-1})\,x_k^{|\lambda|-|\delta|-|\mu|+k-1-|\nu/\mu|}\,,\\
\end{eqnarray*}
where the sum is over all $\nu$ that can be obtained from $\mu$ by adding a vertical strip of size and height at most $k-1$. Now $|\nu/\mu| = |\nu|-|\mu|$ and $|\delta|-|\delta'| = k-1$, so we get
\begin{eqnarray}
\widetilde{\chi}_{\displaystyle {}_{\lambda}}(x_1,\ldots,x_k) & = & \sum_{\mu\,\lhd\,\lambda-\delta}\sum_{\nu}\widetilde{\chi}_{\displaystyle {}_{\nu+\delta'}}(x_1,\ldots,x_{k-1})\,x_k^{|\lambda|-|\nu+\delta'|}\\
& = & \sum_{\nu\in\Lambda^{\!+}_S\,:\,\nu\,\lhd\,\lambda}2^{\nabla(\lambda,\nu)}\,\widetilde{\chi}_{\displaystyle {}_{\nu}}\,x_k^{|\lambda|-|\nu|}\,.
\end{eqnarray}
Iterating as before yields
\begin{equation}
\widetilde{\chi}_{\displaystyle {}_{\lambda}} = \sum_{\lambda^{(1)}\,\lhd\,\cdots\,\lhd\,\lambda^{(k)}\,=\,\lambda}\hspace{-3mm}2^{\nabla(\lambda^{(k)},\lambda^{(k-1)})+\cdots+\nabla(\lambda^{(2)},\lambda^{(1)})}\,\widetilde{\chi}_{\displaystyle {}_{\lambda^{(1)}}}\,x_1^{|\lambda^{(1)}|}\,x_2^{|\lambda^{(2)}|-|\lambda^{(1)}|}\,\cdots\,x_k^{|\lambda^{(k)}|-|\lambda^{(k-1)}|}\,.
\end{equation}
Therefore $\widetilde{V}_\mathcal{D}$ lies in the weight space with weight
\begin{equation}
\left(|\lambda^{(1)}|, |\lambda^{(2)}|-|\lambda^{(1)}|, \ldots, |\lambda^{(k)}|-|\lambda^{(k-1)}|\right)\,,
\end{equation}
which is the same as what we found in equations~\eqref{eqn:twistedGTweight0} and~\eqref{eqn:twistedGTweight1}.
\end{remark}

\section{A chamber complex for the $q$-analogue}

We will assume in this section that we are working in type $A_n$. We will let $\Phi_+ = \{\alpha_1,\ldots,\alpha_N\}$, with $N={n+1\choose 2}$. We will denote by $\mathcal{C}_n$ the chamber complex for the Kostant partition function. For positive integer $k$, we will use the notation $[k]$ for the set $\{1,2,\ldots,k\}$.

\begin{definition}
Let $M$ be a $d \times n$ matrix over the integers, such that $\mathrm{ker}M \cap \mathbb{R}^n_{\geq 0} = 0$. The \emph{vector partition function} (or simply \emph{partition function}) associated to $M$ is the function
\begin{displaymath}
\begin{array}{rccl}
\phi_M : & \mathbb{Z}^d & \longrightarrow & \mathbb{N}\\
& b & \mapsto & |\{x \in \mathbb{N}^n \ : \  Mx = b\}|
\end{array}
\end{displaymath}
\end{definition}

The condition $\mathrm{ker}M \cap \mathbb{R}^n_{\geq 0} = 0$ forces the set $\{x \in \mathbb{N}^n \ : \  Mx = b\}$ to have finite size, or equivalently, the set $\{x \in \mathbb{R}^n_{\geq 0} \ : \  Mx = b\} $ to be compact, in which case it is a polytope $P_b$, and the partition function is the number of integral points (lattice points) inside it.

Also, if we let $M_1, \ldots, M_n$ denote the columns of $M$ (as column-vectors), and $x=(x_1, \ldots, x_n) \in \mathbb{R}^n_{\geq 0}$, then $Mx = x_1M_1 + x_2M_2 + \cdots + x_nM_n$ and for this to be equal to $b$, $b$ has to lie in the cone $\mathrm{pos}(M)$ spanned by the vectors $M_i$. So $\phi_M$ vanishes outside of $\mathrm{pos}(M)$. 

It is well-known that partition functions are piecewise quasipolynomial, and that the domains of quasipolynomiality form a complex of convex polyhedral cones, called the \emph{chamber complex}. Sturmfels gives a very clear explanation in \cite{Sturmfels} of this phenomenon. The explicit description of the chamber complex is due to Alekseevskaya, Gel'fand and Zelevinski$\breve{\i}$ \cite{AGZ}. There is a special class of matrices for which partition functions take a much simpler form. Call an integer $d\times n$ matrix $M$ of full rank $d$ \emph{unimodular} if every nonsingular $d\times d$ submatrix has determinant $\pm 1$. For unimodular matrices, the chamber complex determines domains of polynomiality instead of quasipolynomiality \cite{Sturmfels}.

It is useful for what follows to describe how to obtain the chamber complex of a partition function. Let $M$ be a $d\times n$ integer matrix of full rank $d$ and $\phi_M$ its associated partition function. For any subset $\sigma \subseteq \{1,\ldots,n\}$, denote by $M_\sigma$ the submatrix of $M$ with column set $\sigma$, and let $\tau_\sigma = \mathrm{pos}(M_\sigma)$, the cone spanned by the columns of $M_\sigma$. Define the set $\mathcal{B}$ of \emph{bases} of $M$ to be
\begin{displaymath}
\mathcal{B} = \{\sigma \subseteq \{1,\ldots,n\}\ : \ |\sigma| = d \ \ \textrm{and } \ \mathrm{rank}(M_\sigma) = d\}\,.
\end{displaymath}
$\mathcal{B}$ indexes the invertible $d\times d$ submatrices of $M$. The \emph{chamber complex} of $\phi_M$ is the common refinement of all the cones $\tau_\sigma$, as $\sigma$ ranges over $\mathcal{B}$ (see \cite{AGZ}). A theorem of Sturmfels \cite{Sturmfels} describes exactly how partition functions are quasipolynomial over the chambers of that complex. 

If we let $M_{A_n}$ be the matrix whose columns are the positive roots $\Phi_+^{(A_n)}$ of $A_n$, written in the basis of simple roots, then we can write Kostant's partition function in the matrix form defined above as
\begin{displaymath}
K_{A_n}(v) = \phi_{M_{A_n}}(v)\,.
\end{displaymath}

The following lemma is a well-known fact about $M_{A_n}$ and can be deduced from general results on matrices with columns of $0$'s and $1$'s where the $1$'s come in a consecutive block (see \cite{Schrijver}).

\begin{lemma}
\label{lemma:unimodularity}
The matrix $M_{A_n}$ is unimodular for all $n$. 
\end{lemma}

$M_{A_n}$ unimodular means that the Kostant partition functions for $A_n$ is polynomial instead of quasipolynomial on the cells of the chamber complex. In general, for $M$ unimodular, the polynomial pieces have degree at most the number of columns of the matrix minus its rank (see \cite{Sturmfels}). In our case, $M_{A_n}$ has rank $n$ and as many columns as $A_n$ has positive roots, ${n+1\choose 2}$. Hence the Kostant partition function for $A_n$ is piecewise polynomial of degree at most ${n+1\choose 2}-n = {n\choose 2}$.

We can now state the main result of this section.

\begin{theorem}
The $q$-analogue $K_q(\mu)$ is given by polynomials of degree ${n\choose 2}$ with coefficients in $\mathbb{Q}[q]$ of degree ${n+1\choose 2}$ over the relative interior of the cells of $\mathcal{C}_n$.
\end{theorem}

\begin{proof}
We start at the level of generating functions by observing that
\begin{eqnarray}
\prod_{i=1}^N\frac{1+(q-1)e^{\alpha_i}}{1-e^{\alpha_i}} & = & \underbrace{\left(\prod_{i=1}^N\frac{1}{1-e^{\alpha_i}}\right)}_{\displaystyle\sum_{\mu}K(\mu)e^\mu}\cdot\underbrace{\left(\prod_{i=1}^N(1+(q-1)e^{\alpha_i})\right)}_{\displaystyle\sum_{I\subseteq [N]}(q-1)^{|I|}e^{\alpha_I}}\nonumber\\[2mm]
& = & \sum_{I\subseteq [N]}(q-1)^{|I|}\sum_{\mu}K(\mu)e^{\mu+\alpha_I}\,,
\end{eqnarray}
where $\alpha_I = \displaystyle\sum_{i\in I}\alpha_i$.

Extracting the coefficient of $e^\mu$ in the previous equation gives
\begin{equation}
K_q(\mu) = \sum_{I\subseteq [N]}(q-1)^{|I|}K(\mu-\alpha_I)\,.
\end{equation}

Now,
\begin{eqnarray*}
K(\mu-\alpha_I) & = & \Big|\Big\{{(k_i)}_{i\in [N]}\ :\ \mu-\alpha_I = \sum_{i=1}^N k_i\alpha_i\Big\}\Big|\\
& = & \Big|\Big\{{(k_i)}_{i\in [N]}\ :\ \mu = \sum_{i=1}^N k_i\alpha_i+\sum_{i\in I}\alpha_i\Big\}\Big|\\
& = & \Big|\Big\{{(k_i)}_{i\in [N]}\ :\ \mu = \sum_{i=1}^N k_i\alpha_i, \ \textrm{$k_i\geq 1$ if $i\in I$}\Big\}\Big|\,.
\end{eqnarray*}

Applying inclusion-exclusion, we get
\begin{eqnarray*}
K(\mu-\alpha_I) & = & \Big|\Big\{{(k_i)}_{i\in [N]}\ :\ \mu = \sum_{i=1}^N k_i\alpha_i\Big\}\Big| \quad-\quad \sum_{j\in I}\Big|\Big\{{(k_i)}_{i\in [N]}\ :\ \mu = \sum_{i=1}^N k_i\alpha_i, \ k_j=0\Big\}\Big|\\
& & + \sum_{j_1,j_2\in I, j_1\neq j_2}\Big|\Big\{{(k_i)}_{i\in [N]}\ :\ \mu = \sum_{i=1}^N k_i\alpha_i, \ k_{j_1}=k_{j_2}=0\Big\}\Big| \quad-\quad \ldots \\[5mm]
& = & \sum_{J\subseteq I}(-1)^{|J|}\,\Big|\underbrace{\Big\{{(k_i)}_{i\in [N]}\ :\ \mu = \sum_{i=1}^N k_i\alpha_i, \ \textrm{$k_j=0$ if $j\in J$}\Big\}}_{\displaystyle\Big\{{(k_i)}_{i\in [N]\backslash J}\ :\ \mu = \sum_{i\in [N]\backslash J} k_i\alpha_i\Big\}}\Big|\,.
\end{eqnarray*}

Note that the function
\begin{equation}
\mu \quad\longmapsto\quad \Big|\Big\{{(k_j)}_{j\in J}\ :\ \mu = \sum_{j\in J}k_j\alpha_j\Big\}\Big|
\end{equation}
is a vector partition function, corresponding to the matrix $M_J$ with $\{\alpha_j\,:\,j\in J\}$ as columns. We will denote this function by $K_J(\mu)$. With this notation, we can write
\begin{equation}
\label{eqn:q-analogue}
K_q(\mu) = \sum_{I\subseteq [N]}(q-1)^{|I|}\sum_{J\subseteq I}(-1)^{|J|}K_{[N]\backslash J}(\mu)\,.
\end{equation}

Denote by $\mathcal{C}_J$ the chamber complex associated to the partition function $K_J$. If $M_J$ has full rank $n$, then the base cones whose common refinement is $\mathcal{C}_J$ are the positive hulls of the columns of the nonsingular $n\times n$ submatrices of $M_J$. As these submatrices are also nonsingular $n\times n$ submatrices of $M_{A_n}=M_{[N]}$, the base cones of $\mathcal{C}_J$ are also base cones of $\mathcal{C}_n$, and $\mathcal{C}_J$ is therefore a coarsening of the chamber complex for the Kostant partition function. If $M_J$ does not have full rank $n$, then the base cones of $\mathcal{C}_J$ are the positive hulls of the columns of the maximal-rank submatrices of $M_J$. These cones are faces of the base cones for $M_{A_n}=M_{[N]}$. The complex $\mathcal{C}_J$ will therefore be a coarsening of the restriction of $\mathcal{C}_n$ to the positive hull of the columns of $M_J$ (a lower dimensional complex).   

In view of all this, if $\mu$ stays strictly within any given cell of $\mathcal{C}_n$, it also stays in the same cell of $\mathcal{C}_J$ (for any $J$), and $K_J(\mu)$ is obtained by evaluating the same polynomial attached to that cell of $\mathcal{C}_J$. Hence $K_q(\mu)$ is given by polynomials over the relative interior of the cells of $\mathcal{C}_n$. Since $K_J(\mu)$ has degree at most ${n\choose 2}$ in $\mu$ (see remarks after Lemma~\ref{lemma:unimodularity}), equation~\eqref{eqn:q-analogue} gives that $K_q(\mu)$ is polynomial of degree at most ${n\choose 2}$ in $\mu$ with coefficients of degree at most $N={n+1\choose 2}$ in $q$.
\end{proof}

\subsection{An example: $A_2$}

The chamber complex $\mathcal{C}_2$ has two top-dimensional cones:
\begin{eqnarray*}
\tau_1 & = & \{a_1\alpha_1+a_2\alpha_2\ : \ \textrm{$a_1,a_2>0$ and $a_1>a_2$}\}\,,\\
\tau_2 & = & \{a_1\alpha_1+a_2\alpha_2\ : \ \textrm{$a_1,a_2>0$ and $a_1<a_2$}\}\,,
\end{eqnarray*}
three 1-dimensional cones:
\begin{eqnarray*}
\tau_3 & = & \{a(\alpha_1+\alpha_2)\ : \ a>0\}\,,\\
\tau_4 & = & \{a_1\alpha_1\ : \ a_1>0\}\,,\\
\tau_5 & = & \{a_2\alpha_2\ : \ a_2>0\}\,,\\
\end{eqnarray*}
and the 0-dimensional cone
\begin{displaymath}
\tau_6 = \{0\}\,.
\end{displaymath}

For $\mu=(\mu_1,\mu_2,\mu_3)$ in the root lattice (in particular, $\mu_1+\mu_2+\mu_3=0$), we get
\begin{equation}
K_q(\mu) = \left\{\begin{array}{lcl}
(\mu_1+\mu_2-1)q^3+2q^2 & & \textrm{if $\mu\in\tau_1$}\,,\\[2mm]
(\mu_1-1)q^3 + 2q^2 & & \textrm{if $\mu\in\tau_2$}\,,\\[2mm]
(\mu_1-1)q^3+q^2+q & & \textrm{if $\mu\in\tau_3$}\,,\\[2mm]
q & & \textrm{if $\mu\in\tau_4$ or $\mu\in\tau_5$}\,,\\[2mm]
1 & & \textrm{if $\mu\in\tau_6$}\,,\\[2mm]
0 & & \textrm{otherwise}\,.
\end{array}\right.
\end{equation}

\bigskip

\begin{center}
{\large \textbf{Acknowledgements}}
\end{center}

We would like to thank Richard Stanley for suggesting that the Schur function approach might work, rather than our more complicated approach in terms of Hall-Littlewood polynomials, and also for the observation that the tensor product of two twisted representations can be written as a positive sum (rather than as a virtual sum) of twisted representations. We would also like to thank Shlomo Sternberg and David Vogan for useful discussions and comments.

\medskip

\bibliographystyle{plain}

\bibliography{SignatureQuantization}


\end{document}